\documentclass[12pt,leqno]{article} 
\usepackage{amssymb} 
\usepackage{amsbsy}
\usepackage{amsmath} 
\usepackage{amsfonts} 
\usepackage{amscd}
\usepackage{graphicx,psfrag} 
\usepackage{color} 
 
\newcommand{\email}[1]{{\small E-mail: {\textsf {#1}}}}

\newtheorem{thm}{Theorem} 
\newtheorem{prop}[thm]{Proposition} 
\newtheorem{lem}[thm]{Lemma} 
 
\newtheorem{rem}[thm]{Remark} 
 
\newcommand{\beqn}{\begin{equation}} 
\newcommand{\eeqn}{\end{equation}} 
\newcommand{\bear}{\begin{eqnarray}} 
\newcommand{\eear}{\end{eqnarray}} 
\newcommand{\bean}{\begin{eqnarray*}} 
\newcommand{\eean}{\end{eqnarray*}} 
\newcommand\qed{\hfill $\square$} 
\newcommand{\RR}{\mathbb R}

\newcommand{\e}{\varepsilon} 
\begin{document} 

\title{Finite time blow-up for a one-dimensional\\ 
quasilinear parabolic-parabolic chemotaxis system}

\author{Tomasz Cie\'slak\footnote{Institute of Applied Mathematics, Warsaw University, Banacha 2, 02-097 Warszawa, Poland. \email{T.Cieslak@impan.gov.pl}} \kern8pt \& \kern8pt
Philippe Lauren\c cot\footnote{Institut de 
Math\'ematiques de Toulouse, CNRS UMR~5219, Universit\'e de Toulouse, 118 route de Narbonne, F--31062 Toulouse Cedex 9, 
France. \email{Philippe.Laurencot@math.univ-toulouse.fr}}}
\date{\today}
\maketitle

\begin{abstract}
Finite time blow-up is shown to occur for solutions to a one-dimensional quasilinear parabolic-parabolic chemotaxis system as soon as the mean value of the initial condition exceeds some threshold value. The proof combines a novel identity of virial type with the boundedness from below of the Liapunov functional associated to the system, the latter being peculiar to the one-dimensional setting. 
\end{abstract}

%%%%%%%%%%%%%%%%%%%%%%%%%%%%%%%%%%
\section{Introduction}\label{i}
%%%%%%%%%%%%%%%%%%%%%%%%%%%%%%%%%%

%\setcounter{thm}{0}
%\setcounter{equation}{0}

We study the possible occurrence of blow-up in finite time for solutions to a one-dimensional parabolic system modeling chemotaxis \cite{KS70}. More precisely, we consider the Keller-Segel chemotaxis model with nonlinear diffusion which describes the space and time evolution of a population of cells moving under the combined effects of diffusion (random motion) and a directed motion in the direction of high gradients of a chemical substance (chemoattractant) secreted by themselves. If $u\ge 0$ and $v$ denote the density of cells and the (rescaled) concentration of chemoattractant, respectively, the Keller-Segel model with nonlinear diffusion reads
\bear
\label{a1}
\partial_t u = {\rm div}\left( a(u)\ \nabla u - u\ \nabla v \right) & \;\;\mbox{ in }\;\; & (0,\infty)\times\Omega\,, \\
\label{a2}
\e\ \partial_t v = D\, \Delta v - \gamma\ v + u - M & \;\;\mbox{ in }\;\; & (0,\infty)\times\Omega\,, \\
\label{a3}
a(u)\ \partial_\nu u = \partial_\nu v = 0 & \;\;\mbox{ on }\;\; & (0,\infty)\times\partial\Omega\,, \\
\label{a4}
(u,v)(0) = (u_0,v_0) & \;\;\mbox{ in }\;\; & \Omega\,.
\eear
In general, $\Omega$ is an open bounded subset of $\RR^N$, $N\ge 1$, with smooth boundary $\partial\Omega$, $a$ is a smooth non-negative function, and the parameters $\e$, $D$, $\gamma$, and $M$ are non-negative real numbers with $D>0$ and $M>0$. In addition, the initial data $u_0$ and $v_0$ satisfy
\beqn
\label{a5}
u_0\ge 0\,, \;\; \int_\Omega u_0(x)\ dx = |\Omega|\ M\,, \;\;\mbox{ and }\;\; \int_\Omega v_0(x)\ dx = 0\,.
\eeqn
The constraints \eqref{a5} ensure in particular that a solution $(u,v)$ to \eqref{a1}-\eqref{a4} satisfies (at least formally) the same properties for positive times, that is, 
\beqn
\label{a6}
u(t)\ge 0\,, \;\; \int_\Omega u(t,x)\ dx = |\Omega|\ M\,, \;\;\mbox{ and }\;\; \int_\Omega v(t,x)\ dx = 0\,.
\eeqn

The main feature of \eqref{a1} is that it involves a competition between the diffusive term ${\rm div}\left( a(u)\ \nabla u \right)$ (spreading the population of cells) and the chemotactic drift term $- {\rm div}\left( u\ \nabla v \right)$ (concentrating the population of cells) that may lead to the blow-up in finite time of the solution to \eqref{a1}-\eqref{a4}. The possible occurrence of such a singular phenomenon is actually an important mathematical issue in the study of \eqref{a1}-\eqref{a4} which is also relevant from a biological point of view: indeed, it corresponds to the experimentally observed concentration of cells in a narrow region of the space which is a preamble to a change of state of the cells. From a mathematical point of view, the blow-up issue has been the subject of several studies in the last twenty years, see the survey \cite{Ho03} and the references therein. 

Still, it is far from being fully understood, in particular when $\e>0$ (the so-called parabolic-parabolic Keller-Segel model). In that case, the only finite time blow-up result available seems to be that of Herrero \& Vel\'azquez who showed in \cite{HV96b,HV97} that, when $\Omega$ is a ball in $\RR^2$, $D=1$, and $a\equiv1$, there are $M>8\pi$ and radially symmetric solutions $(u,v)$ to \eqref{a1}-\eqref{a4} which blow up in finite time. These solutions are constructed as small perturbations of time rescaled stationary solutions to \eqref{a1}-\eqref{a4} and a similar result is also true when $\e=0$ \cite{HV96a}. The result  in \cite{HV97} actually goes far beyond the mere occurrence of blow-up in finite time as the shape of the blow-up profile is also identified. Recall that the condition $M>8\pi$ is necessary for the finite blow-up to take place: indeed, it is shown in \cite{NSY97} that, if $\Omega$ is a ball in $\RR^2$, $D=1$, and $a\equiv1$, radially symmetric solutions to \eqref{a1}-\eqref{a4} are global as soon as $M<8\pi$. We refer to \cite{GZ98,NSY97} for additional global existence results when $\Omega $ is a bounded domain in $\RR^2$, $\e>0$, and $a\equiv1$. In \cite{Ho02} the existence of unbounded solutions is shown for $\e>0$ and $a\equiv 1$, but it is not known whether the blow-up takes place in finite or infinite time. The same approach is employed in \cite{hw} to obtain unbounded solutions to quasilinear Keller-Segel systems, still without knowing whether the blow-up time is finite or infinite. The finite time blow-up result proved in this paper (Theorem~\ref{tha1}) is thus the first one of this kind for quasilinear parabolic-parabolic Keller-Segel systems. 

In contrast, for the parabolic-elliptic Keller-Segel system corresponding to $\e=0$, several finite time blow-up results are available. There is thus a discrepancy between the two cases $\e>0$ and $\e=0$ which may be explained as follows. On the one hand, as observed in \cite{JL92} when $\e=0$, $\Omega$ is a ball of $\RR^2$, $a\equiv1$, and $u_0$ is radially symmetric, it is possible to reduce \eqref{a1}-\eqref{a4} to a single parabolic equation for the cumulative distribution function 
$$
U(t,r) := \int_{B(0,r)} u(t,x)\ dx\,.
$$
Finite time blow-up is then shown with the comparison principle by constructing appropriate subsolutions. This approach was extended to nonlinear diffusions (non-constant $a$) and arbitrary space dimension $N\ge 1$ in \cite{CW08}. On the other hand, it has been noticed in \cite{BiN94,Na95} that, still for $a\equiv1$, the moment $M_k$ of $u$ defined by 
$$
M_k(t) := \int_\Omega |x|^k\ u(t,x)\ dx\,, \qquad k\in (0,\infty)\,,
$$
satisfies a differential inequality which cannot hold true for all times for a suitably chosen value of $k>0$, for it would imply that $u$ reaches negative values in finite time in contradiction with \eqref{a6}. In contrast to the previous approach, this is an obstructive method which provides no information on the blow-up profile and is somehow reminiscent of the celebrated virial identity available for the nonlinear Schr\"odinger equation (see, e.g., \cite[Section~6.5]{Ca03} and the references therein). Nevertheless, it applies to more general sets $\Omega$ \cite{Na00,Na01,SS01}. We recently develop further this technique to establish finite time blow-up of radially symmetric solutions to \eqref{a1}-\eqref{a4} with $\e=0$ in a ball of $\RR^N$, $N\ge 2$, when the diffusion is nonlinear \cite{CLxx}, the main idea being to replace the moments by nonlinear functions of the cumulative distribution function $U$. For a related model in $\RR^N$ with nonlinear diffusion $a(u)=m\ u^{m-1}$, $m>1$, finite time blow-up  results were recently established in  \cite{BCLxx,Su06} by looking at the evolution of the second moment $M_2$. 

\medskip

Coming back to the parabolic-parabolic Keller-Segel system \eqref{a1}-\eqref{a4} ($\e>0$), it seems unlikely that the first approach described above (reduction to a single equation) could work and the purpose of this paper is to show that finite time blow-up results can be established by the second approach in the one-dimensional case ($N=1$). More precisely, we consider the initial-boundary value problem 
\bear
\label{a7}
\partial_t u = \partial_x\left( a(u)\ \partial_x u - u\ \partial_x v \right) & \;\;\mbox{ in }\;\; & (0,\infty)\times (0,1)\,, \\
\label{a8}
\e\ \partial_t v = D\, \partial_x^2 v - \gamma\ v + u - M & \;\;\mbox{ in }\;\; & (0,\infty)\times (0,1)\,, \\
\label{a9}
a(u)\ \partial_x u = \partial_x v = 0 & \;\;\mbox{ on }\;\; & (0,\infty)\times\{0,1\}\,, \\
\label{a10}
(u,v)(0) = (u_0,v_0) & \;\;\mbox{ in }\;\; & (0,1)\,,
\eear
and assume that 
\beqn
\label{a11}
\e>0\,, \quad D>0\,, \quad \gamma\ge 0\,, \quad M>0\,,
\eeqn
and the initial data $(u_0,v_0)\in W^{1,2}(0,1;\RR^2)$ satisfy
\beqn
\label{a12}
u_0\ge 0\,, \;\; \int_0^1 u_0(x)\ dx = M\,, \;\;\mbox{ and }\;\; \int_0^1 v_0(x)\ dx = 0\,.
\eeqn
We further assume that $a\in\mathcal{C}^2(\RR)$ and that there are $p\in (1,2]$, and $c_1>0$ such that 
\beqn
\label{a13}
0<a(r) \le c_1\ (1+r)^{-p}\;\;\mbox{ for }\;\; r\ge 0\,.
\eeqn
Our main result then reads as follows.

\begin{thm}\label{tha1}
Assume that the parameters $\e$, $D$, $\gamma$, $M$, the initial data $(u_0,v_0)$, and the function $a$ fulfil the conditions \eqref{a11}, \eqref{a12}, and \eqref{a13}, respectively. Then there is a unique classical maximal solution 
$$
(u,v)\in \mathcal{C}([0,T_m)\times [0,1];\RR^2)\cap \mathcal{C}^{1,2}((0,T_m)\times [0,1];\RR^2)
$$
to \eqref{a7}-\eqref{a10} with maximal existence time $T_m\in (0,\infty]$. It also satisfies
\beqn
\label{a14}
u(t,x)\ge 0\,, \;\; \int_0^1 u(t,x)\ dx = M\,, \;\;\mbox{ and }\;\; \int_0^1 v(t,x)\ dx = 0
\eeqn
for $(t,x)\in [0,T_m)\times [0,1]$. Introducing 
\bear
\nonumber
F(z_1,z_2) & := & c_1\ \left( 1+M \right) + \frac{M^2}{2D}  + z_1 + M\ z_2 + \frac{D+\gamma}{2}\ z_2^2\,, \\
\label{a15}
\mathcal{P}_q(z_1,z_2,z_3) & := & \left( 1 + \frac{\gamma}{D} + \frac{\gamma}{M}\ z_2 + \frac{M^{q-2}}{4qD}\ z_3 \right)\ F(z_1,z_2) \\
\nonumber
& + & \frac{c_1 (q-1) q^{(q-2)/q} D}{(p-1) M^{p-1}}\ F(z_1,z_2)^{(q-2)/q} - \frac{M^q}{q(q+1)}
\eear
and
$$
m_q(0) := \frac{1}{q}\ \int_0^1 \left( \int_0^x u_0(y)\ dy \right)^q\ dx\,,
$$
for $(z_1,z_2,z_3)\in [0,\infty)^3$ and $q\geq 2$, we have $T_m<\infty$ as soon as $\mathcal{P}_q\left( m_q(0),\|v_0\|_{H^1},\e M \right)<0$ for some finite $q\in (2,2/(2-p)]$. In particular, if $u_0$ is such that 
\beqn
\label{a16}
\mathcal{P}_q\left( m_q(0),0,0 \right)<0 \;\;\mbox{ for some finite }\;\; q\in (2,2/(2-p)]\,,
\eeqn 
there is $\vartheta>0$ such that $\e M\in (0,\vartheta)$ and $\|v_0\|_{H^1}<\vartheta$ imply that $\mathcal{P}_q\left( m_q(0),\|v_0\|_{H^1},\e M \right)<0$ and thus $T_m<\infty$. 
\end{thm}

There are functions $u_0$ satisfying \eqref{a12} and \eqref{a16} if $M$ is sufficiently large. Indeed, observe that 
\bean
\mathcal{P}_q\left( 0,0,0 \right) & = & \left( 1 + \frac{\gamma}{D} \right)\ \left( c_1\ \left( 1+M \right) + \frac{M^2}{2D} \right) \\
\nonumber
& + & \frac{c_1 (q-1) q^{(q-2)/q} D}{(p-1) M^{p-1}}\ \left( c_1\ \left( 1+M \right) + \frac{M^2}{2D} \right)^{(q-2)/q} - \frac{M^q}{q(q+1)}
\eean
is negative for sufficiently large $M$ as $q>2$. Given such an $M>0$ and choosing the function $u_0(x)= 2M\ \max{\{ x+\delta-1 , 0 \}}/\delta^2$, $x\in (0,1)$, we have $m_q(0)=(2M)^q \delta/(2q+1)$ and $\mathcal{P}_q\left( m_q(0),0,0 \right)<0$ for $\delta>0$ small enough. 
In fact, if $u_0$ fulfils \eqref{a16}, then the same computation as the one leading to Theorem~\ref{tha1} shows that the corresponding solution to the parabolic-elliptic Keller-Segel system ($\e=0$) blows up in a finite time and the last assertion of Theorem~\ref{tha1} states that this property remains true for the parabolic-parabolic Keller-Segel system $(\e>0)$ provided $\e$ and $v_0$ are small, that is, in a kind of neighbourhood of the parabolic-elliptic case.

\begin{rem}\label{rea2}
The growth condition required on $a$ in \eqref{a13} is seemingly optimal: indeed, it is proved in \cite{CW08} that $T_m=\infty$ if $a(r)\ge c_0\ (1+r)^{-p}$ for some $p<1$ and $\e=0$, and the proof is likely to extend to the case $\e>0$. 
\end{rem}

The proof of Theorem~\ref{tha1} relies on two properties of the Keller-Segel system \eqref{a7}-\eqref{a10}: first, there is a Liapunov functional \cite{GZ98} which is bounded from below in the one-dimensional case \cite{CW08} and which provides information on the time derivative of $v$. This will be the content of Section~\ref{wplf} where we also sketch the proof of the local well-posedness of \eqref{a7}-\eqref{a10}. We next derive an identity of virial type for the $L^q$-norm of the indefinite integral of $u$ in Section~\ref{ftbu} which involves in particular the time derivative of $v$. The information obtained on this quantity in the previous section then allow us to derive a differential inequality for the $L^q$-norm of the indefinite integral of $u$ for a suitable value of $q$ which cannot be satisfied for all times if the parameters $\e$, $D$, $\gamma$, $M$, and the initial data $(u_0,v_0)$ are suitably chosen.

%%%%%%%%%%%%%%%%%%%%%%%%%%%%%%%%%%
\section{Well-posedness and Liapunov functional}\label{wplf}
%%%%%%%%%%%%%%%%%%%%%%%%%%%%%%%%%%

%\setcounter{thm}{0}
%\setcounter{equation}{0}

In this section, we establish the local well-posedness of \eqref{a7}-\eqref{a10} in $W^{1,2}(0,1;\RR^2)$ and recall the availability of a Liapunov functional for this system \cite{GZ98}. To this end, we assume that 
\beqn
\label{b1}
0<a\in\mathcal{C}^2(\RR)
\eeqn
and define $b\in \mathcal{C}^2((0,\infty))$ by
\beqn
\label{b2}
b(1)=b'(1):=0 \;\;\mbox{ and }\;\; b''(r) := \frac{a(r)}{r} \;\;\mbox{ for }\;\; r>0\,.
\eeqn

\begin{prop}\label{prb1}
Assume that the parameters $\e$, $D$, $\gamma$, $M$, and the function $a$ fulfil \eqref{a11} and \eqref{b1}, respectively. Given initial data $(u_0,v_0)\in W^{1,2}(0,1;\RR^2)$ satisfying \eqref{a12}, there is a unique classical maximal solution 
$$
(u,v)\in \mathcal{C}([0,T_m)\times [0,1];\RR^2)\cap \mathcal{C}^{1,2}((0,T_m)\times [0,1];\RR^2)
$$
to \eqref{a7}-\eqref{a10} with maximal existence time $T_m\in (0,\infty]$ and $(u,v)$ satisfies \eqref{a14} for $t\in [0,T_m)$. In addition, if $T_m<\infty$, we have
\beqn
\label{b5}
\lim_{t\to T_m} \left( \|u(t)\|_\infty + \|v(t)\|_\infty \right) = \infty\,.
\eeqn
\end{prop}

\noindent\textbf{Proof.} We define $\tilde{a}\in \mathcal{C}^2(\RR^2;\mathcal{M}_2(\RR))$ by 
$$
\tilde{a}(y) = \left( \tilde{a}^{m,n}(y) \right)_{1\le m,n\le 2} := \left( 
\begin{array}{cc}
D & 0 \\
 & \\
-y_2 & a(y_2) 
\end{array}
\right) \;\;\mbox{ for }\;\; y=(y_1,y_2)\in\RR^2
$$
and introduce the operators  
\bean
\mathcal{A}(y) z & := & - \partial_x \left( \tilde{a}(y)\ \partial_x z \right) \,,\\
\left( \mathcal{B}(y) z(0),\mathcal{B}(y) z(1) \right) & := & \left( - \tilde{a}(y)\ \partial_x z(0) , \tilde{a}(y)\ \partial_x z(1) \right)\,,
\eean
and the function
$$
f(y) := \left( 
\begin{array}{c}
-\gamma\ y_1 + y_2 - M \\
 \\
0
\end{array}
\right)
$$
with $z=(z_1,z_2)$. With this notation, an abstract formulation of \eqref{a7}-\eqref{a10} reads
\bean
\partial_t z + \mathcal{A}(z) z & = & f(z)\,,\\
\mathcal{B}(z) z & = & 0\,, \\
z(0) & = & (v_0,u_0)\,,
\eean
with $z=(v,u)$ and we aim at applying the theory developed in \cite{Am93}. Owing to \eqref{a11} and \eqref{b1}, $\tilde{a}(y)$ is a positive definite matrix for all $y\in\RR^2$ and we infer from \cite[Section~4]{Am93} that the boundary-value operator $(\mathcal{A},\mathcal{B})$ is normally elliptic. It then follows from \cite[Theorems~14.4 \&~14.6]{Am93} that \eqref{a7}-\eqref{a10} has a unique classical maximal solution 
$$
(v,u):=z\in \mathcal{C}([0,T_m)\times [0,1];\RR^2)\cap \mathcal{C}^{1,2}((0,T_m)\times [0,1];\RR^2)
$$
for some $T_m\in (0,\infty]$. In addition, $\tilde{a}^{2,1}(y_1,0)=0$ for $y_1\in\RR$ and we deduce from \cite[Theorem~15.1]{Am93} that $u(t,x)\ge 0$ for $(t,x)\in [0,T_m)\times [0,1]$. The property \eqref{a14} then readily follows from \eqref{a7}-\eqref{a10} and \eqref{a12} by integration. As for the last statement \eqref{b5}, it is a consequence of the lower triangular structure of the matrix $\tilde{a}$ and \cite[Theorem~15.5]{Am93}. \qed

\medskip

We next proceed as in \cite{GZ98} to check the availability of a Liapunov functional for \eqref{a7}-\eqref{a10}. 

\begin{lem}\label{leb2}
Assume that the parameters $\e$, $D$, $\gamma$, $M$, and the function $a$ fulfil \eqref{a11} and \eqref{b1}, respectively. Given initial data $(u_0,v_0)\in W^{1,2}(0,1;\RR^2)$ satisfying \eqref{a12} and such that $b(u_0)\in L^1(0,1)$, the corresponding classical solution $(u,v)$ to \eqref{a7}-\eqref{a10} satisfies
\beqn
\label{b6}
L(u(t),v(t)) + \e\ \int_0^t \|\partial_t v(s)\|_2^2\ ds \le L(u_0,v_0) \;\;\mbox{ for }\;\; t\in [0,T_m)\,,
\eeqn
where
\beqn
\label{b7}
L(u,v):= \int_0^1 \left( b(u) - u v + \frac{D}{2}\ |\partial_x v|^2 + \frac{\gamma}{2}\ |v|^2 \right)\ dx\,.
\eeqn
\end{lem}

\noindent\textbf{Proof.} It follows from \eqref{a7}-\eqref{a9} that
\bear
\nonumber
\frac{d}{dt} L(u,v) & = & \int_0^1 \left( b'(u) - v \right)\ \partial_t u\ dx + \int_0^1 \left( D\ \partial_x v\ \partial_x\partial_t v + (\gamma\ v - u)\ \partial_t v \right)\ dx \\
\nonumber
& = & - \int_0^1 \left( b''(u)\ \partial_x u - \partial_x v \right)\ \left( a(u)\ \partial_x u - u\ \partial_x v \right)\ dx \\
\nonumber
& + & \int_0^1 \partial_t v\ \left( -D\ \partial_x^2 v + \gamma\ v - u \right)\ dx  \\
\nonumber
& = & - \int_0^1 u\ \left| \partial_x \left( b'(u)-v \right) \right|^2\ dx - \int_0^1 \left( M + \e\ \partial_t v \right)\ \partial_t v\ dx \\
\label{b6b}
& \le & - \e\ \|\partial_t v\|_2^2\,,
\eear
the last inequality being a consequence of \eqref{a14}. Integrating the previous inequality with respect to time gives \eqref{b6}. \qed

\medskip

We next take advantage of the one-dimensional setting to show that $L$ is bounded from below without prescribing growth conditions on $a$. This fact has already been observed in \cite{CW08} and is peculiar to the one-dimensional case. Indeed, as shown in \cite{GZ98,Ho02}, the occurrence of blow-up is closely related to the unboundeness of the Liapunov functional. 

\begin{lem}\label{leb3}
Assume that the parameters $\e$, $D$, $\gamma$, $M$, and the function $a$ fulfil \eqref{a11} and \eqref{b1}, respectively. Given initial data $(u_0,v_0)\in W^{1,2}(0,1;\RR^2)$ satisfying \eqref{a12} and such that $b(u_0)\in L^1(0,1)$, the corresponding classical solution $(u,v)$ to \eqref{a7}-\eqref{a10} satisfies
\beqn
\label{b8}
L(u(t),v(t)) \ge - \frac{M^2}{2D} \;\;\mbox{ for }\;\; t\in [0,T_m)\,.
\eeqn
\end{lem}

\noindent\textbf{Proof.} Owing to \eqref{a14}, the Poincar\'e inequality ensures that $\|v(t)\|_\infty\le \|\partial_x v(t)\|_2$ for $t\in [0,T_m)$ so that
$$
\int_0^1 u(t)\ v(t)\ dx \le \|v(t)\|_\infty\ \| u(t)\|_1 \le \|\partial_x v(t)\|_2\ \|u(t)\|_1\,.
$$
We use again \eqref{a14} as well as the non-negativity of $b$ to conclude that 
$$
L(u(t),v(t)) \ge \frac{D}{2}\ \|\partial_x v(t)\|_2^2 - M\ \|\partial_x v(t)\|_2 = \frac{D}{2}\ \left( \|\partial_x v(t)\|_2 - \frac{M}{D} \right)^2 - \frac{M^2}{2D}
$$
for $t\in [0,T_m)$, from which \eqref{b8} readily follows. \qed

%%%%%%%%%%%%%%%%%%%%%%%%%%%%%%%%%%
\section{Finite time blow-up}\label{ftbu}
%%%%%%%%%%%%%%%%%%%%%%%%%%%%%%%%%%

%\setcounter{thm}{0}
%\setcounter{equation}{0}

As already mentioned, the main novelty in this paper is a new identity of virial type which is the cornerstone of the proof that blow-up takes place in finite time under suitable assumptions. Specifically, we assume that the parameters $\e$, $D$, $\gamma$, $M$, and the function $a$ fulfil the conditions \eqref{a11} and \eqref{a13}, respectively. Recalling the definition \eqref{b2} of $b$, we deduce from \eqref{a13} that  
\beqn
\label{ap1}
b(r) \le c_1\ \left( r\ \ln{r} - r + 1 \right)\ \mathbf{1}_{[0,1]}(r) + \frac{c_1(r-1)}{p}\ \mathbf{1}_{[1,\infty)}(r) \le c_1\ (1+r)\,, \quad r\ge 0\,.
\eeqn
We also define 
\beqn
\label{ap2}
A(r) := - \int_r^\infty a(s)\ ds\,, \quad r\ge 0\,,
\eeqn
and infer from \eqref{a13} that $A$ is well-defined and satisfies
\beqn
\label{ap3}
0\le - A(r)\ r \le \frac{c_1}{p-1}\ r^{2-p}\,, \quad r\ge 0\,.
\eeqn

Consider next initial data $(u_0,v_0)\in W^{1,2}(0,1;\RR^2)$ satisfying \eqref{a12}. If $(u,v)$ denotes the corresponding classical solution to \eqref{a7}-\eqref{a10} given by Proposition~\ref{prb1}, we define the cumulative distribution functions $U$ and $V$ by
\beqn
\label{c1}
U(t,x) := \int_0^x u(t,y)\ dy \;\;\mbox{ and }\;\; V(t,x) := \int_0^x v(t,y)\ dy\eeqn
for $(t,x)\in [0,T_m)\times [0,1]$. It readily follows from \eqref{a7}-\eqref{a9} and \eqref{a14} that $(U,V)$ solves
\bear
\label{c2}
\partial_t U = \partial_x A(u) - u \ \partial_x v & \;\;\mbox{ in } & (0,T_m)\times (0,1)\,, \\
\label{c3}
\e\ \partial_t V = D\ \partial_x v - \gamma\ V + U - Mx & \;\;\mbox{ in } & (0,T_m)\times (0,1)\,, 
\eear
the function $A$ being defined in \eqref{ap2}, and
\beqn
\label{c4}
U(t,0)=M-U(t,1)=0 \;\;\mbox{ and }\;\; V(t,0)=V(t,1)=0\,, \quad t\in [0,T_m)\,.
\eeqn

\begin{lem}\label{lec1z}
Introducing $m_q(t):=\|U(t)\|_q^q/q$ for $q\ge 2$, we have
\bear
\label{c5z}
\frac{dm_q}{dt} & = & \frac{M}{D}\ m_q - \frac{M^{q+1}}{q(q+1)D} + M^{q-1}\ A(u(t,1)) - (q-1)\ \int_0^1 U^{q-2}\ u\ A(u)\ dx \\
\nonumber
& + & \frac{\e}{qD}\ \int_0^1 U^q\ \partial_t v\ dx - \frac{\gamma}{D}\ \int_0^1 U^{q-1}\ u\ V\ dx
\eear
for $t\in [0,T_m)$.
\end{lem}

\noindent\textbf{Proof.} We infer from \eqref{c2}, \eqref{c3}, and \eqref{c4} that 
\bean
\frac{dm_q}{dt} & = & \left[ U^{q-1}\ A(u) \right]_{x=0}^{x=1} - (q-1)\ \int_0^1 U^{q-2}\ u\ A(u)\ dx \\
& - & \frac{1}{D}\ \int_0^1 u\ U^{q-1}\ \left( \e\ \partial_t V + \gamma\ V - U + Mx \right)\ dx \\
& = & M^{q-1}\ A(u(t,1)) - (q-1)\ \int_0^1 U^{q-2}\ u\ A(u)\ dx - \frac{\e}{qD}\ \left[ U^q\ \partial_t V \right]_{x=0}^{x=1} \\
& + & \frac{\e}{qD}\ \int_0^1 U^q\ \partial_t v\ dx - \frac{\gamma}{D}\ \int_0^1 U^{q-1}\ u\ V\ dx + \frac{1}{(q+1)D}\ \left[ U^{q+1} \right]_{x=0}^{x=1} \\
& - & \frac{M}{qD}\ \left[ U^q\ x \right]_{x=0}^{x=1} + \frac{M}{D}\ m_q\\
& = & M^{q-1}\ A(u(t,1)) - (q-1)\ \int_0^1 U^{q-2}\ u\ A(u)\ dx + \frac{\e}{qD}\ \int_0^1 U^q\ \partial_t v\ dx \\
& - &  \frac{\gamma}{D}\ \int_0^1 U^{q-1}\ u\ V\ dx - \frac{M^{q+1}}{q(q+1)D} + \frac{M}{D}\ m_q\,,
\eean
which is the expected identity. \qed

\medskip

At this point, we notice that the solution to the ordinary differential equation $D\ \dot{X} = M\ X - (M^{q+1}/(q(q+1)))$ (obtained by neglecting several terms in \eqref{c5z}) is given by 
$$
X(t) = \frac{M^q}{q(q+1)} + e^{Mt/D}\ \left( X(0) - \frac{M^q}{q(q+1)} \right)\,,
$$
and thus vanishes at a finite time if $X(0)<M^q/(q(q+1))$. If a similar argument could be used for \eqref{c5z}, we would obtain a positive time $t_0$ such that $m_q(t_0)=0$ which clearly contradicts the properties of $U(t_0)$: indeed, by \eqref{c1} and \eqref{c4}, $x\mapsto U(t_0,x)$ is continuous with $U(t_0,1)=M$. Consequently, the solution $(u,v)$ to \eqref{a7}-\eqref{a10} no longer exists at this time $t_0$ and blow-up shall have occurred at an earlier time, thus establishing Theorem~\ref{tha1}. For this approach to work, we shall of course control the other terms on the right-hand side of \eqref{c5z} which will in turn give rise to the blow-up criterion stated in Theorem~\ref{tha1}. The latter is actually a simple consequence of the following result:

\begin{thm}\label{thc2z}
Assume that the parameters $\e$, $D$, $\gamma$, $M$, and the initial data $(u_0,v_0)$ are such that 
\beqn
\label{cfinalz}
E\left( m_q(0) + L(u_0,v_0) + \frac{M^2}{2D} \right)  < 0
\eeqn
for some finite $q\in (2,2/(2-p)]$, where
$$
E(z):= \left( 1 + \frac{\gamma}{D} + \frac{\gamma}{M}\ \|v_0\|_{H^1} + \frac{\e M^{q-1}}{4qD}\right)\ z + \frac{c_1 (q-1) q^{(q-2)/q} D}{(p-1)M^{p-1}}\ z^{(q-2)/q} - \frac{M^q}{q(q+1)}
$$ 
for $z\ge 0$. Then $T_m<\infty$.
\end{thm}

\noindent\textbf{Proof.} The starting point of the proof being the identity \eqref{c5z}, we first derive upper bounds for the terms on the right-hand side of \eqref{c5z} involving $A$, $\e$, and $\gamma$. Thanks to \eqref{ap3} and the non-negativity of $U$, it follows from the H\"older inequality that
\bean
& & M^{q-1}\ A(u(t,1)) - (q-1)\ \int_0^1 U^{q-2}\ u\ A(u)\ dx \\ 
& \le & \frac{c_1 (q-1)}{p-1}\ \int_0^1 U^{q-2}\ u^{2-p}\ dx \\
& \le & \frac{c_1 (q-1) q^{(q-2)/q}}{(p-1)}\ m_q^{(q-2)/q}\ \left( \int_0^1 u^{((2-p)q)/2}\ dx \right)^{2/q}\,.
\eean
Since $q\in  (2,2/(2-p)]$, we may use the Jensen inequality and \eqref{a14} to conclude that
\beqn
\label{c6z}
M^{q-1}\ A(u(t,1)) - (q-1)\ \int_0^1 U^{q-2}\ u\ A(u)\ dx \le \frac{c_1 (q-1) q^{(q-2)/q}}{(p-1)}\ M^{2-p}\ m_q^{(q-2)/q}\,.
\eeqn
Next, to estimate the term involving $\gamma$, we adapt an argument from \cite{Na95} and first claim that 
\beqn
\label{c7}
V(t,x)\ge V_m(t,x) := \frac{M}{6D}\ (x^3-x) + h(t,x)\,, \quad (t,x)\in [0,T_m)\times [0,1]\,,
\eeqn
where $h$ denotes the unique solution to 
\bear
\label{c8}
& & \e\ \partial_t h - D\ \partial_x^2 h + \gamma\ h = 0\,, \quad (t,x)\in  (0,\infty)\times (0,1)\,,\\
\label{c9}
& & h(t,0)=h(t,1)=0\,, \quad t\in (0,\infty)\,,\\
& & h(0,x) = \min{\left\{ V(0,x) + \frac{M}{6D}\ (x-x^3) , 0 \right\}} \le 0 \,, \quad x\in (0,1)\,.
\eear
Indeed, $V_m\le V$ on $[0,T_m)\times\{0,1\}$ and $\{0\}\times [0,1]$, and it follows from the non-negativity of $U$ and the negativity of $h$ that
\bean
\e\ \partial_t V_m - D\ \partial_x^2 V_m + \gamma\ V_m & = & \e\ \partial_t h - M x - D\ \partial_x^2 h + \frac{M \gamma}{6D}\ (x^3-x) + \gamma\ h \\
& \le & -Mx \le U - Mx = \e\ \partial_t V - D\ \partial_x^2 V + \gamma\ V\,.
\eean
The comparison principle then implies \eqref{c7}. We next infer from \eqref{c7} and the non-negativity of $u$ and $U$ that 
\bean
-  \frac{\gamma}{D}\ \int_0^1 U^{q-1}\ u\ V\ dx & \le & -  \frac{\gamma}{D}\ \int_0^1 U^{q-1}\ u\ V_m\ dx \\
& = & -  \frac{\gamma}{qD}\ \left[ U^q\ V_m \right]_{x=0}^{x=1} +   \frac{\gamma}{qD}\ \int_0^1 U^q\ \partial_x V_m\ dx\\
& \le & \frac{\gamma}{D}\ \left( \frac{M}{2D} + \|\partial_x h\|_\infty \right)\ m_q\,.
\eean
We next note that $\partial_x h$ also solves \eqref{c8} with homogeneous Neumann boundary conditions, the latter property being a consequence of \eqref{c8} and \eqref{c9}. Since 
$$
\left| \partial_x h(0,x) \right| \le \left| v_0(x) + \frac{M}{6D}\ (1-3x^2) \right| \le \|v_0\|_\infty + \frac{M}{3D}\,,
$$
the comparison principle and the non-negativity of $\gamma$ warrant that $\|\partial_x h(t)\|_\infty \le \|v_0\|_\infty + (M/3D)$ for $t\ge 0$. Consequently, recalling the Sobolev embedding $\|v_0\|_\infty\le \|v_0\|_{H^1}$, we end up with 
\beqn
\label{c11z}
-  \frac{\gamma}{D}\ \int_0^1 U^{q-1}\ u\ V\ dx \le \frac{\gamma M}{D^2}\ \left( 1 + \frac{D}{M}\ \|v_0\|_{H^1} \right)\ m_q\,.
\eeqn
We finally infer from \eqref{a14}, \eqref{c1}, \eqref{c4}, and the H\"older inequality that
\beqn
\label{c12z}
\frac{\e}{qD}\ \int_0^1 U^q\ \partial_t v\ dx \le \frac{\e M^{q/2}}{qD}\ \int_0^1 U^{q/2}\ |\partial_t v|\ dx \le \frac{\e M^{q/2}}{q^{1/2}D}\ m_q^{1/2}\ \|\partial_t v\|_2\,.
\eeqn

It now follows from \eqref{c5z}, \eqref{c6z}, \eqref{c11z}, and \eqref{c12z} that 
\bean
\frac{dm_q}{dt} & \le & \frac{M}{D}\ \left[ \left( 1 + \frac{\gamma}{D} + \frac{\gamma}{M}\ \|v_0\|_{H^1} \right)\ m_q + \frac{c_1 (q-1) q^{(q-2)/q} D}{(p-1) M^{p-1}}\ m_q^{(q-2)/q} - \frac{M^q}{q(q+1)} \right] \\
& + & \frac{\e M^{q/2}}{q^{1/2}D}\ m_q^{1/2}\ \|\partial_t v\|_2 \\
& \le & \frac{M}{D}\ E(m_q) - \frac{\e M^q}{4qD^2}\ m_q + \frac{\e M^{q/2}}{q^{1/2}D}\ m_q^{1/2}\ \|\partial_t v\|_2\,.
\eean
Owing to \eqref{a12} and \eqref{ap1}, we have $b(u_0)\in L^1(0,1)$ and it follows from \eqref{b6b}, \eqref{b8}, and the above inequality that
\bean
\frac{d}{dt}\left( m_q + L(u,v) + \frac{M^2}{2D} \right) & \le & \frac{M}{D}\ E(m_q) - \frac{\e M^q}{4qD^2}\ m_q + \frac{\e M^{q/2}}{q^{1/2}D}\ m_q^{1/2}\ \|\partial_t v\|_2 - \e\ \|\partial_t v\|_2^2 \\
& = & \frac{M}{D}\ E(m_q) - \e\ \left( \|\partial_t v\|_2 - \frac{M^{q/2}}{2 q^{1/2} D}\ m_q^{1/2} \right)^2 \\
& \le & \frac{M}{D}\ E(m_q)\,.
\eean
Using now the monotonicity of $E$ and \eqref{b8}, we end up with
$$
\frac{d}{dt}\left( m_q + L(u,v) + \frac{M^2}{2D} \right) \le  \frac{M}{D}\ E\left( m_q + L(u,v) + \frac{M^2}{2D} \right)\,.
$$

Assume now for contradiction that $T_m=\infty$. The previous inequality and \eqref{cfinalz} then warrant that there is a time $t_0>0$ such that $m_q(t_0)+L(u(t_0),v(t_0))+(M^2/2D)=0$ and hence $m_q(t_0)=0$ by \eqref{b8}. This in turn implies that $U(t_0,x)=0$ for all $x\in [0,1]$ and contradicts \eqref{c4}. Consequently, $T_m<\infty$. \qed

\medskip

The remaining step towards Theorem~\ref{tha1} is to use the properties of $a$ to simplify the condition \eqref{cfinalz} derived in Theorem~\ref{thc2z}. 

\noindent\textbf{Proof of Theorem~\ref{tha1}.}
It follows from \eqref{a12}, \eqref{ap1}, and the Sobolev embedding  $\|v_0\|_\infty\le \|v_0\|_{H^1}$ that
\bean
L(u_0,v_0) + \frac{M^2}{2D} & \le & \int_0^1 \left( c_1\ (1+u_0) + \frac{D}{2}\ |\partial_x v_0|^2 + \frac{\gamma}{2}\ |v_0|^2 + u_0\ \|v_0\|_\infty \right)\ dx + \frac{M^2}{2D}\\
& \le & c_1\ (1+M) + \frac{M^2}{2D} + \frac{D+\gamma}{2}\ \|v_0\|_{H^1}^2 + M\ \|v_0\|_{H^1} \\
& = & F\left( m_q(0), \|v_0\|_{H^1} \right) - m_q(0)\,,
\eean
the function $F$ being defined in Theorem~\ref{tha1}. Therefore, 
$$
E\left( m_q(0) + L(u_0,v_0) + \frac{M^2}{2D} \right) \le (E \circ F)\left( m_q(0), \|v_0\|_{H^1} \right) = \mathcal{P}_q\left( m_q(0), \|v_0\|_{H^1}, \e M \right)\,,
$$
and the condition $\mathcal{P}_q\left( m_q(0),\|v_0\|_{H^1},\e M\right)<0$ clearly implies \eqref{cfinalz} and hence $T_m<\infty$. \qed

\medskip

\textbf{Acknowledgement.} This paper was prepared during T.~Cie\'slak's one-month visit at the Institut de Math\'ematiques de Toulouse, Universit\'e Paul Sabatier. T. Cie\'slak would like to express his gratitude for the invitation, support, and hospitality.

%%%%%%%%%%%%%%%%%%%%%%%%%%%%%%%%%%

\end{document}